\documentclass[preprint,final,1p]{elsarticle}

\usepackage{amsmath,amssymb}
\usepackage{epsfig,color,hyperref,setspace}

\journal{Journal of Power Sources}

\usepackage{filecontents}

\begin{document}

\begin{frontmatter}

\title{Study of voltage cycling conditions on Pt oxidation and dissolution 
in polymer electrolyte fuel cells}

\author[mymainaryaddress,mysecondaryaddress]{V.A.~Kovtunenko}
\ead{victor.kovtunenko@uni-graz.at}

\author[mythirdaddress]{L.~Karpenko-Jereb\corref{mycorrespondingauthor}}
\ead{karpenkojereb@gmail.com}
\cortext[mycorrespondingauthor]{Corresponding author}

\address[mymainaryaddress]{Institute for Mathematics and Scientific Computing,
Karl-Franzens University of Graz, NAWI Graz, Heinrichstr. 36, 8010 Graz, Austria}
\address[mysecondaryaddress]{Lavrentyev Institute of Hydrodynamics,
Siberian Division of the Russian Academy of Sciences, 630090 Novosibirsk, Russia}
\address[mythirdaddress]{Institute of Electronic Sensor Systems, 
Graz University of Technology, Inffeldgasse 10/II, Graz 8010, Austria}

\begin{abstract}
This paper is devoted to study the electrochemical behavior of Pt catalyst 
in a polymer electrolyte fuel cell at various operating conditions and 
at different electric potential difference (also known as voltage) cycling 
applied in accelerated stress tests. 
The degradation of platinum is considered with respect to 
the Pt ion dissolution and the Pt oxide coverage of catalyst 
described by a one-dimensional model. 
In the model, degradation rate increases with temperature and 
decreasing particle diameter of Pt nano-particles. 
The theoretical study of the underlying diffusion system with the nonlinear 
reactions is presented by analytical methods and 
gives explicit solutions through a first integral of the ODE system. 
Numerical tests are obtained using a second order implicit-explicit scheme. 
The computer simulation shows that the lifetime of the catalyst depends 
on the voltage profile and the upper potential level. 
By this Pt mass loss is more significant at the membrane surface than at the gas diffusion layer. 
\end{abstract}


\end{frontmatter}

\section{Introduction}\label{sec1}

At present, the polymer electrolyte fuel cells are very extensively developed 
as power sources for application in portable computers, heating systems 
as well for leisure yachts, aircraft and vehicles \cite{EK/17,KJA/18,SS/18,Wee/07}. 
In commercial PEM fuel cells, the most used material for catalyst is platinum. 
Effective usage of this valuable metal is challenging task in elaboration 
of the modern fuel cell systems with high durability and long lifetime. 

Mathematical modelling helps to better understand phenomena causing 
the chemical degradation on the Pt surface and allows predicting a decline 
in electro-chemical surface area (ECSA) of the catalyst and effectivity 
of the fuel cell \cite{KJSFT/16,KB/11,ZYLZSWZ/20}. 
In the last decades, numerous publications have been devoted to develop mathematical 
models prediction Pt oxidation and dissolution in the catalyst layer of PEMFCs. 
In 2007, Zhang et al \cite{ZLGLG/07} suggested a simple model accounting 
Pt dissolution and further precipitation within the membrane of PEMFCs. 
One year late, Bi and Fuller \cite{BF/08} published a dynamic model 
applied to estimate Pt mass loss during potential cycling. 
In 2009, Holby et al. \cite{HSSHM/09} developed a kinetic model 
taking into account effect of Pt particle size on the rate of the degradation. 
It was demonstated that, due to rapid changes in the Gibbs--Thomson energy, 
particle size effects dominate degradation for 2 nm particles 
but play almost no role for 5 nm particles. 
In 2015 Hiraoka et al. \cite{HKMM/15} proposed a model for the Pt particle 
growth based on the Gibbs--Thomson equation and simulated change 
in the Pt particle size distribution during electric potential cycling. 
Using the developed model, the authors investigated an effect of 
high potential limit and Pt particle diameter on particle size growth. 
The results showed that lowing high potential in the cycling 
the Pt particle size grows more slowly. 
A decrease in particle size accelerates Pt dissolution. 
In the same year, Li et al. \cite{LMGAW/15} published a mathematical model 
taking into account effect of RH, T and catalyst layer thickness on the loss 
in electrochemical active area of the Pt catalyst. 
The simulation results were in good agreement with experiments. 
The conducted simulation for a triangle potential cycle demonstrated that 
thinning the cathode catalyst layer would induce more rapid ECSA loss. 
ECSA increases with rising temperature and higher relative humidity. 
In 2018, using a simplified model Baricci et al. \cite{BBYGMC/18} studied 
changes in roughness through the catalyst layer thickness. 
Non-uniform degradation is observed in the catalyst layer consequently 
to the formation of a platinum depleted region next to the membrane, 
which, according to the model, results from diffusion and precipitation 
of dissolved platinum into the membrane. 
Recently, Koltsiva et al. \cite{KVSFB/18} suggested a novel model, 
which for the first time considers five phenomena simultaneously 
proceeding on the Pt/C catalyst surface: 
platinum nanoparticles electrochemical dissolution, 
particle growth due to Ostwald ripening, 
migration of nanoparticles along the carbon support, 
coalescence of fine particles, 
diffusion of platinum ions in the ionomer. 

The present paper is focused on study the effect of kind of voltage cycling, where 
voltage is commonly adopted as electric potential difference versus a reference of 0 (V),  
on the electro-chemical behavior of Pt catalyst in a polymer electrolyte fuel cell (PEFC). 
For this purpose we have developed a mathematical model based 
on the physico-chemical data of the Pt dissolution and oxidation reactions 
and taking account diffusion of building Pt ions into the ionomer membrane. 
The two degradation phenomena of the platinum ion (Pt$^{2+}$) 
dissolution and the platinum oxide (PtO) formation of Pt catalyst layer (coverage) 
in polymer-electrolyte membrane fuel cells (PEMFC) 
are studied theoretically at different voltage cycling conditions. 
To describe these phenomena, the degradation model due to 
Holby \cite{HM/12,HSSHM/09} is utilized, 
which is one-dimensional (1d) across the catalyst layer (CL) thickness 
and accounts for diffusion of Pt ions. 
For mathematical modeling of interface reactions in multi-phase media 
within complete electrokinetic Poisson--Nernst--Planck equations, 
see e.g. \cite{FK/15,GGK/18,KZ/18}. 
For mathematical approaches which are suitable to describe and to test 
a mechanical degradation due to fracture phenomena, 
we refer to \cite{HKK/07,IKR/20,KK/00}. 

The paper is organized in the following way: 
after the introduction we describe theoretical approach developed, where  
the part ``Degradation model of Pt catalyst'' displays geometrical, physical and 
chemical properties of the catalyst layer and the ionomer membrane used for further simulation; 
``Theoretical and numerical methods'' introduces the mathematical and the numerical models; 
``Results'' presents simulation results and their discussion; 
``Discussion'' and ``Conclusion'' summarize the most important findings of the current study. 
The mathematical model is given by a coupled system of nonlinear 
reaction-diffusion equations with modified Butler--Volmer reaction rates 
for three unknown variables: 
Pt$^{2+}$ concentration, Pt particle diameter, and PtO coverage ratio. 
For its numerical solution we use a second order implicit-explicit scheme 
following \cite{ARS/97}. 
Neglecting diffusion of Pt ions, 
in \ref{A} the resulting nonlinear reaction equations are reduced 
to two unknown variables with the help of a first integral of the system, 
while in \ref{B} an example of analytical solution is constructed. 

\section{Degradation model of Pt catalyst}\label{sec2}

\paragraph{Approximations} 

The developed model of the Platinum on carbon (Pt/C) degradation 
is based on the following approximations: 

\begin{itemize}
\item[1.]
This is a one-dimensional and a dynamic model. 
\item[2.]
The model considers two layers of a polymer electrolyte fuel cell (PEFC): 
the cathode catalyst layer (CL) with length $L_{\rm CL}$ and 
the polymer electrolyte membrane (PEM) with length $L_{\rm PEM}$. 
The catalyst layer is filled with catalyst particles: spherical Pt nano-particles 
placed on C-support bound with the membrane by perfluorinated sulfonated ionomer, 
the membrane is made from the same ionomer. 
\item[3.]  
The degradation of Pt nano-particles are caused by Pt oxidation and Pt dissolution, 
which are described by the following electro-chemical reactions: 
\begin{subequations}\label{1}
\begin{equation}\label{1a}
{\rm Pt}_{\rm (s)}\longleftrightarrow {\rm Pt}^{2+}_{\rm (aq)} +2{\rm e}^-,
\end{equation}
\begin{equation}\label{1b}
{\rm Pt}_{\rm (s)} +{\rm H}_2{\rm O}_{\rm (aq)}\longleftrightarrow 
{\rm PtO}_{\rm (s)} +2{\rm H}^+_{\rm (aq)} +2{\rm e}^-.
\end{equation}
\end{subequations}
\item[4.]  
The rates of the degradation reactions \eqref{1} are simulated by 
a modified Butler--Volmer  equation (see \eqref{3} and \eqref{4}) 
taking into account Gibbs--Thomson's effect: 
dependence of a surface potential on nano-particle size as well as 
influence of the surface potential on the potential gradient in the system. 
\item[5.]  
The platinum ions occurring due to the Pt dissolution \eqref{1} diffuse through 
the ionomer phase of the catalyst layer into the polymer electrolyte membrane. 
The diffusion of Pt ions into the gas diffusion layer is impossible 
because this layer does not possess ionic conductivity. 
This determines boundary conditions \eqref{5e} for the Pt ions diffusion \eqref{5a}.  
\item[6.]  
On balance, the model takes into account an effect of potential gradient, 
Pt particle size,  temperature,  relative humidity as well as other phenomena 
(see parameter gathered in Table~\ref{tab1})
on the degradation rates of platinum catalyst and allows calculating 
the platinum ion concentration (${\rm Pt}^{2+}$), the particle diameter, 
and the platinum oxide (PtO) coverage ratio from governing relations \eqref{5}. 
\end{itemize}

With respect to the relative humidity in point 6 we remark the following. 
In our model we consider an effect of pH on the reaction rate of Pt oxidation. 
The pH depends on the dissociation degree of sulfonyl groups of the ionomer. 
In the simulated cases we suppose that the sulfonyl groups 
are completely dissociated and $pH=0$. 
Generally, pH is function of proton concentration, in PEMFCs pH depends on relative humidity. 

\paragraph{Model variables} 

\begin{figure}[hbt!]
\begin{center}
\epsfig{file=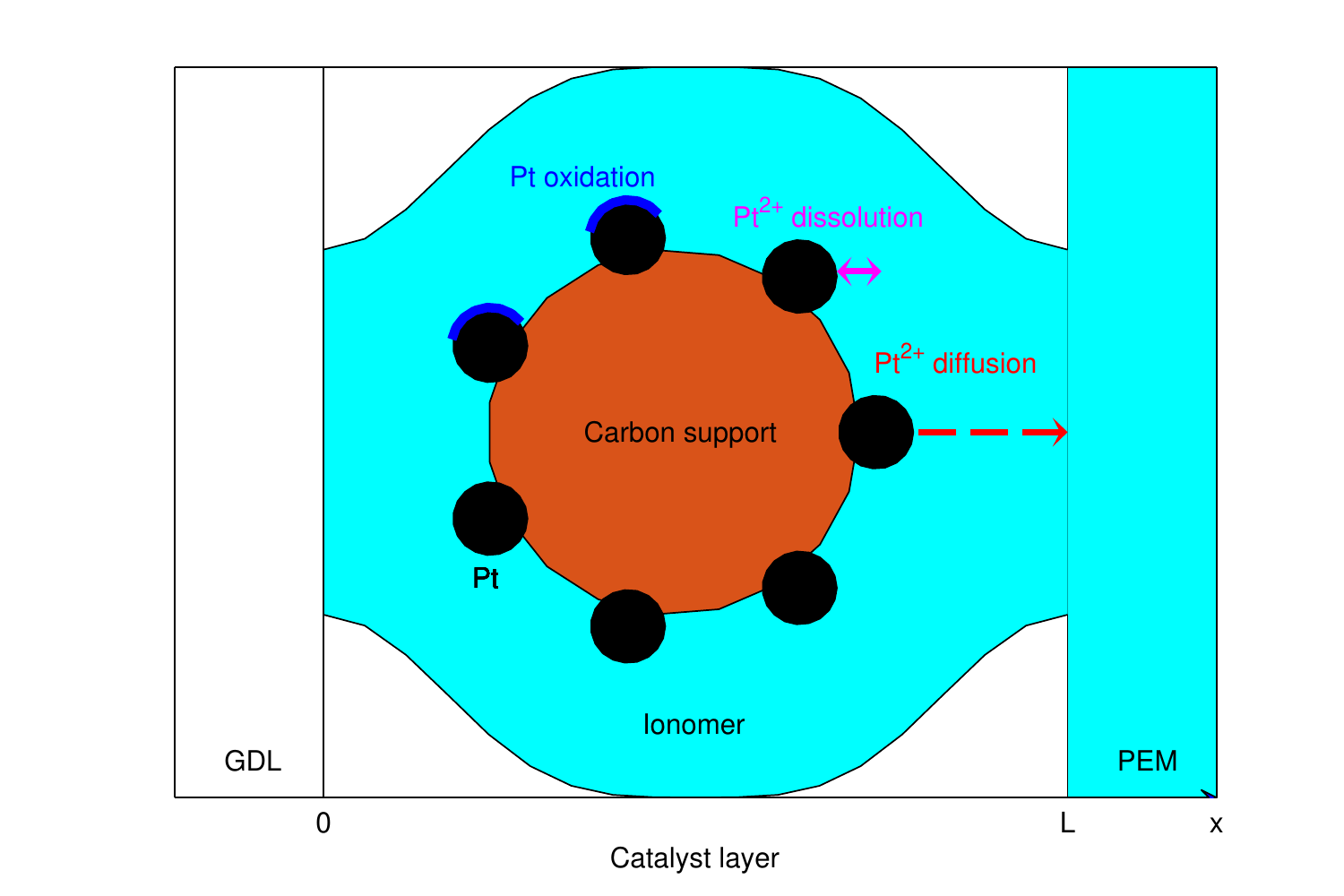,width=.6\textwidth,angle=0}
\caption{The model of degradation.}
\label{degradation}
\end{center}
\end{figure}

For a semi-infinite cathode catalyst layer of thickness $L$, 
we introduce a spatial variable $x\in[0,L]$ across the CL such that 
one end point $x=0$ corresponds to the CL-gas diffusion layer (GDL) interface, 
and the other end point $x=L$ confirms the CL-PEM interface. 
By this, we allow non-steady state operating conditions for CL 
with respect to time $t\ge0$ and the spatial dependence on 
$x\in[0,L]$ due to diffusion phenomena. 
The model of degradation is sketched in Figure~\ref{degradation}. 
We emphasize that Figure~\ref{degradation} shows a schematic of the model configuration, 
when the Pt particles are fully surrounded by ionomer on the carbon support. 
The real scenario at Pt particles has a partial coverage by ionomer (for ion transfer) 
and a partial coverage by carbon network (for electron transfer) 
and a partial open space to air or gas (for oxygen gas diffusion). 
However, those factors are not addressed in the current model. 

Let Pt particles be hemispheres posed with density $\rho_{\rm Pt}$ 
and loading $p_{\rm Pt}$ on a carbon support, 
such that the Pt volume fraction across CL can be estimated as follows 
$\varepsilon_{\rm Pt} =\frac{p_{\rm Pt}/\rho_{\rm Pt}}{L}$ 
with $\varepsilon_{\rm Pt}<1$. 
The Pt nanoparticle is assumed to be spherical of diameter $d_{\rm Pt}$ 
and volume $V_{\rm Pt} =\frac{4}{3} \pi (\frac{d_{\rm Pt}}{2})^3$. 
Then the Pt number concentration in CL 
is $N_{\rm Pt} =\frac{\varepsilon_{\rm Pt}}{V_{\rm Pt}}$. 
For the parameter values from Table~\ref{tab1}, 
$\varepsilon_{\rm Pt}\approx 2$\%, 
$V_{\rm Pt}\approx 1.5\times10^{-20}$ $({\rm cm}^3)$, 
and $N_{\rm Pt}\approx 1.32\times10^{18}$ $(1/{\rm cm}^3)$. 
 
The unknown constituents entering \eqref{1} 
are the ${\rm Pt}^{2+}$ concentration $c$ ({\rm mol/cm$^3$}), 
Pt particle diameter $d$ ({\rm cm}), and PtO coverage ratio $\theta$, 
that are the time-space dependent functions such that 
\begin{equation}\label{2}
c(t,x)\ge0,\quad d(t,x)\ge0,\quad 0\le\theta(t,x)\le1.
\end{equation}
Based on a modified Butler--Volmer equation, 
the following reaction rates in units of mol/(cm$^2\cdot$s)
are established in \cite{HM/12}: 
for the Pt ion dissolution \eqref{1a}: 
\begin{subequations}\label{3}
\begin{equation}\label{3a}
r_{\rm dissol}(c, d, \theta, V) =B_1(d, \theta) e^{(1-\beta_1) 
B_4(d, \theta) V} -c B_2(d, \theta) e^{-\beta_1 B_4(d, \theta) V},
\end{equation}
where the quantities $B_1$ in mol/(cm$^2\cdot$s), 
$B_2$ in cm/s, $B_4$ in C/J, and $\gamma_0$ in {\rm J/cm$^2$} are 
\begin{equation}\label{3b}
\begin{split}
&B_1(d, \theta) =\nu_1 \mathit{\Gamma} (1 -\theta) e^{\frac{1}{R T}  
(-H_{1,{\rm fit}} -n F (1 - \beta_1) (U_{\rm eq} 
-\frac{4 \mathit{\Omega} \gamma_0(\theta)}{n F d}) )},\\ 
&B_2(d, \theta) ={\textstyle\frac{\nu_2 \mathit{\Gamma}}{c_{\rm ref}}} (1 -\theta) 
e^{\frac{1}{R T} (-H_{1,{\rm fit}} +n F \beta_1 (U_{\rm eq} 
-\frac{4 \mathit{\Omega} \gamma_0(\theta)}{n F d}) )},\\ 
&B_4(d, \theta) ={\textstyle\frac{F}{R T}} (n 
-{\textstyle\frac{4\mathit{\Omega} \mathit{\Gamma} n_2 \theta}{d}}),\\
&\gamma_0(\theta) =\gamma +\mathit{\Gamma} R T \bigl( \theta 
\ln( {\textstyle\frac{\nu_2^\star}{\nu_1^\star}} 10^{-2 pH}) 
+\theta {\textstyle\frac{2 n_2 F U_{\rm fit} +\omega \theta}{2 R T}}  
+\theta \ln ({\textstyle\frac{\theta}{2}}) 
+(2 -\theta) \ln (1 -{\textstyle\frac{\theta}{2}}) \bigr);
\end{split}
\end{equation}
\end{subequations}
for the Pt oxide coverage \eqref{1b}: 
\begin{multline}\label{4}
r_{\rm oxide}(\theta, V) =\mathit{\Gamma} e^{-\frac{1}{R T} (H_{2,{\rm fit}} 
+\lambda \theta)} \bigl( \nu_1^\star  (1 -{\textstyle\frac{\theta}{2}}) 
e^{-\frac{n_2 F (1 - \beta_2)}{R T}  (U_{\rm fit} 
+\frac{\omega \theta}{n_2 F} ) +(1-\beta_2) \frac{n_2 F}{R T} V}\\
-\nu_2^\star 10^{-2 pH} e^{\frac{n_2 F \beta_2}{R T}  (U_{\rm fit} 
+\frac{\omega \theta}{n_2 F} ) -\beta_2 \frac{n_2 F}{R T} V} 
\bigr) .
\end{multline}
The terms in \eqref{3} and \eqref{4} are rearranged in such a way to express 
explicitly the dependence of the governing relations on the voltage $V$.

\begin{figure}[hbt!]
\begin{center}
\hspace*{-15mm}
\epsfig{file=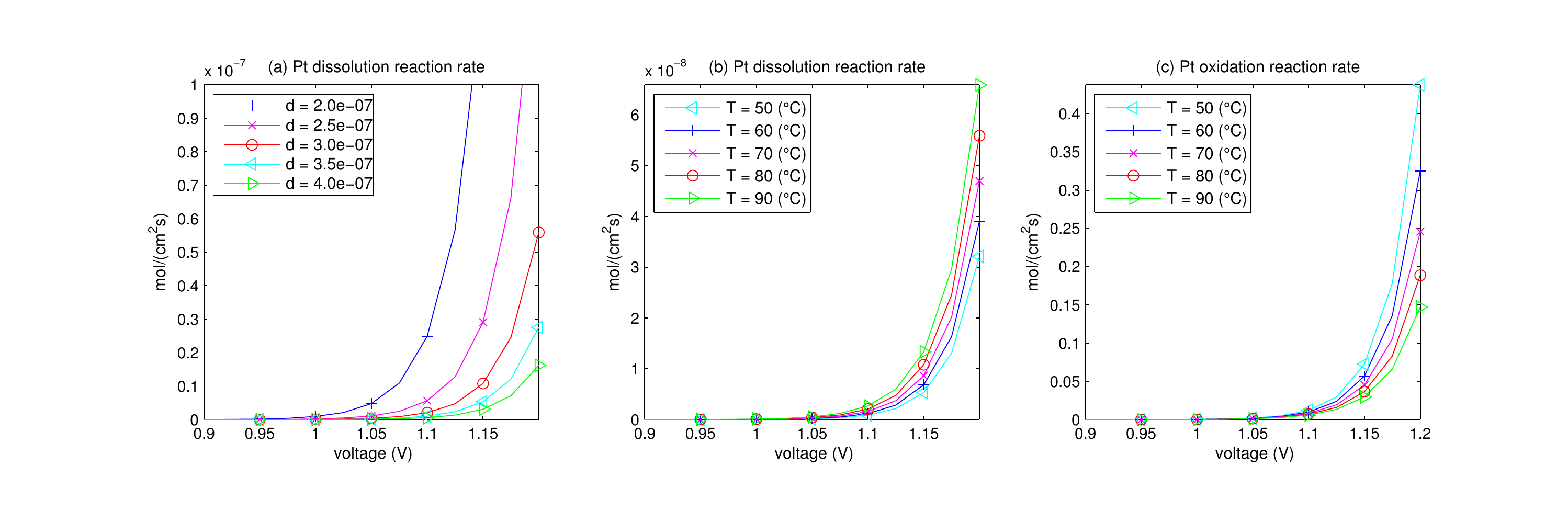,width=1.2\textwidth,angle=0}
\caption{The reaction rates $r_{\rm dissol}$ (a), (b), and 
$r_{\rm oxide}$ (c) for fixed $c$ and $\theta$.}
\label{fig_rateVdT}
\end{center}
\end{figure}

For the parameters from Table~\ref{tab1}, the reactions rates are illustrated 
in Figure~\ref{fig_rateVdT} with respect to varying $V\in[0.9,1.2]$ (V). 
This range corresponds to the operating conditions. 
In Figure~\ref{fig_rateVdT}, there are fixed 
$c=3\times10^{-10}$ ({\rm mol/cm$^3$}) and $\theta=0$, 
the five curves $V\mapsto r_{\rm dissol}(c, d, \theta, \,\cdot\,)$ 
in plot (a) correspond to five Pt particle diameters 
selected equidistantly in the range of $d\in[2,4]\times10^{-7}$ (cm), 
while $r_{\rm oxide}$ is independent of $d$. 
For fixed $d=3\times10^{-7}$ (cm),
the curves $V\mapsto r_{\rm dissol}(c, d, \theta, \,\cdot\,)$ in (b) 
and $V\mapsto r_{\rm oxide}(c, \,\cdot\,)$ in (c) 
are plotted when varying the temperature $T\in[323.15,363.15]$ (K). 
Figure~\ref{fig_rateVdT} demonstrates rates of the forward reactions: 
platinum dissolution (a, b) and formation of platinum oxide (c) 
as function of applied potential gradient at different diameters of Pt particle (a) 
and at different temperature (b, c). 
As seen from the figures the rate of platinum dissolution increases 
with decrease of the diameter of the platinum particles. 
The rate of the both reactions grow with increasing temperature. 

Here it is worth noting that the rates of dissolution and re-deposition reactions are not the same. 
Figure~\ref{fig_rateVdT} shows the rates only forward reactions. 
The model considers dissolution and re-deposition reactions, 
as well as building of platinum oxide (forward reaction, platinum oxidation) 
and reduction of platinum oxide to platinum (backward reaction). 

Since the dependence of $r_{\rm dissol}$ on the Pt concentration $c$ in \eqref{3a} 
is linear, for large values of $c$ the backward dissolution reaction rate 
(having the negative sign) dominates over forward reactions (having the positive sign). 
In these examples, we depict the forward dissolution reaction rates when varying $V$. 
In Figure~\ref{fig_rateVdT} (a) we observe decay of the dissolution 
reaction rates when the particle diameter $d$ increases.  
From Figure~\ref{fig_rateVdT} (b) and (c) we conclude that increasing 
the temperature $T$ follows growth of the dissolution reaction rate 
and decay of the oxidation reaction rate. 

\section{Theoretical and numerical methods}\label{sec3}

\paragraph{Governing relations} 

For a given voltage $V$ in \eqref{3} and \eqref{4}, 
in \cite{LMGAW/15} the following system 
of reaction-diffusion equations is formulated: 
find a triple $(c, d, \theta)$ satisfying \eqref{2} such that 
\begin{subequations}\label{5}
\begin{equation}\label{5a}
{\textstyle\frac{\partial c}{\partial t}} -\sqrt{\varepsilon} D_{\rm Pt} 
{\textstyle\frac{\partial^2 c}{\partial x^2}} =B_3 d^2 r_{\rm dissol}(c, d, \theta)
\quad \text{for $t>0$, $x\in(0,L)$},
\end{equation}
where $B_3 ={\textstyle\frac{\pi N_{\rm Pt}}{2 \varepsilon}}$ 
(1/cm$^3$) is denoted for short, 
\begin{equation}\label{5b}
{\textstyle\frac{\partial d}{\partial t}} =-\mathit{\Omega} 
\,r_{\rm dissol}(c,d,\theta)\quad \text{for $t>0$, $x\in(0,L)$},
\end{equation}
\begin{equation}\label{5c}
{\textstyle\frac{\partial \theta}{\partial t}} 
+{\textstyle\frac{2 \theta}{d}} {\textstyle\frac{\partial d}{\partial t}} 
={\textstyle\frac{r_{\rm oxide}(\theta)}{\mathit{\Gamma}}}
\quad \text{for $t>0$, $x\in(0,L)$};
\end{equation}
which is endowed with the initial conditions:  
\begin{equation}\label{5d}
c =0,\quad d =d_{\rm Pt},\quad \theta =0\quad \text{as $t=0$, $x\in[0,L]$};
\end{equation}
and the mixed Neumann--Dirichlet boundary conditions:  
\begin{equation}\label{5e}
{\textstyle\frac{\partial c}{\partial x}} =0\quad\text{as $t>0$, $x=0$};
\quad c =0\quad \text{as $t>0$, $x=L$}.
\end{equation}
\end{subequations}
The first equality in \eqref{5e} implies no-flux condition at the CL-GDL interface, 
the second condition at the CL-membrane interface assumes 
that the dissolved ${\rm Pt}^{2+}$ concentration goes to zero. 

Neglecting the dependence on the space variable $x$, 
thus omitting the diffusion term (the second one) in the left-hand side 
of \eqref{5a} and the boundary conditions \eqref{5e}, 
the resulting ordinary differential equations (ODE) system 
is studied theoretically in \ref{A}. 
In fact, the problem is reduced to the two unknowns 
by finding a first integral of the system. 
In \ref{B} a particular solution is constructed analytically under specific assumptions. 
The exact solution is used to test numerical solvers for the problem. 
Next we will investigate the initial boundary value problem \eqref{5} numerically. 

\paragraph{Numerical algorithm}

In order to solve the nonlinear reaction-diffusion equation \eqref{5a}, 
below we develop a second order implicit-explicit (IMEX2) scheme 
following \cite{ARS/97}. 

Let the half-strip $[0,\infty)\times[0,L]$ be meshed  with 
temporal points $(t^0,\dots,t^M,\dots)$, where $t_0=0$, and with 
$N+1$ spacial points $(x_0,\dots,x_N)$, where $x_0 =0$ and $x_N=L$. 
For $l=1,\dots,M$, we look for a triple of discrete functions 
\begin{equation}\nonumber
c^l_h =(c^l_0,\dots,c^l_N),\quad d^l_h =(d^l_0,\dots,d^l_N),\quad 
\theta^l_h =(\theta^l_0,\dots,\theta^l_N),
\end{equation} 
with given $c^0_h =\theta^0_h =0$ and $d^0_h =d_{\rm Pt}$
according to the initial condition \eqref{5d}. 
On the space mesh, forward 
$D^+ c^l_i ={\textstyle\frac{c^l_{i+1} -c^l_i}{x_{i+1} -x_i}}$
and backward $D^- c^l_i ={\textstyle\frac{c^l_i -c^l_{i-1}}{x_i -x_{i-1}}}$ 
differences for $i=1,\dots,N-1$ 
are used for the standard approximation of the second-order derivative 
\begin{equation}\nonumber
[D^- D^+] c^l_i ={\textstyle\frac{1}{x_i -x_{i-1}}}
\bigl( {\textstyle\frac{c^l_{i+1} -c^l_i}{x_{i+1} -x_i}} 
-{\textstyle\frac{c^l_i -c^l_{i-1}}{x_i -x_{i-1}}} \bigr).
\end{equation} 

We discretize \eqref{5} and iterate for $l=1,\dots,M$ with the time size 
$\tau^l =t^l -t^{l-1}$ two implicit-explicit equations as follows
\begin{subequations}\label{6}
\begin{equation}\label{6a}
\begin{cases}
c^{l-1/2}_h -w \tau^l \sqrt{\varepsilon} D_{\rm Pt} [D^- D^+] c^{l-1/2}_h
=c^{l-1}_h +w \tau^l B_3 (d^l_h)^2 
r_{\rm dissol}(c^{l-1}_h, d^l_h, \theta^l_h)\\[2ex]
c^l_h =c^{l-1}_h +\tau^l \bigl( 
\sqrt{\varepsilon} D_{\rm Pt} [D^- D^+] c^{l-1/2}_h 
+B_3 (d^l_h)^2 r_{\rm dissol}(c^{l-1/2}_h, d^l_h, \theta^l_h) \bigr),
\end{cases}
\end{equation}
where the IMEX parameter $w=0.5$ is set in the first equation, and  
\begin{equation}\label{6b}
d^l_h =d^{l-1}_h -\tau^l \mathit{\Omega} 
\,r_{\rm dissol}(c^{l-1}_h, d^{l-1}_h, \theta^{l-1}_h),
\end{equation}
\begin{equation}\label{6c}
\theta^l_h =\theta^{l-1}_h +\tau^l \bigl( 
{\textstyle\frac{r_{\rm oxide}(\theta^{l-1}_h)}{\mathit{\Gamma}}} 
+{\textstyle\frac{2 \mathit{\Omega} \theta^{l-1}_h}{d^{l-1}_h}} 
\,r_{\rm dissol}(c^{l-1}_h, d^{l-1}_h, \theta^{l-1}_h) \bigr).
\end{equation}
The diffusion-reaction equations in \eqref{6a} are endowed by 
the boundary conditions according to \eqref{5e}: 
\begin{equation}\label{6d}
c^{l-1/2}_1 =c^{l-1/2}_0,\quad c^l_1 =c^l_0,\quad c^{l-1/2}_N =c^l_N =0.
\end{equation}
\end{subequations}
The standard TDMA algorithm is applied for inversion of a tridiagonal matrix 
in the implicit equation (the first one) in \eqref{6a}. 
For solution of the nonlinear reaction equations \eqref{6b} and \eqref{6c},  
we apply the standard fourth order Runge--Kutta (RK4) method. 

In the numerical examples reported further we set the uniform mesh 
of the time step size $\tau =[10^{-4},10^{-2}]$ (s), 
and the space step size $h =\frac{L}{10}$ (cm) when $N=10$. 
We note that impulse switching of voltage (see 
Figure~\ref{fig_voltage_rate_total} plot (b)) 
requires the smaller time step $\tau = 10^{-4}$ for stable calculation, 
thus increasing the computational complexity. 
The time step $\tau$ was fixed during the iteration. 
The $(\tau, h)$-step choice is conform to the fact, that for large time steps 
the CFL-condition may be violated, thus leading to numerical instabilities. 
The instability appears in such manner that the oxide coverage 
as well as the particle diameter and the Pt2+ concentration becoming negative. 

\section{Results}\label{sec4}

\paragraph{Simulation setup} 

We start with parameter values given in Table~\ref{tab1} at 
the temperature $T = 353.15$ (K) and constant voltage $V = 0.65$ (V) 
that will be used for numerical simulation. 
\begin{table}[hbt!]
{\small
\begin{center}
\begin{tabular}{|l|l|l|p{0.55\textwidth}|l|}\hline
Symbol & Value & Units & Description & Ref.\\\hline
\multicolumn{5}{|l|}{Catalyst layer parameters}\\\hline
$L$ & $1\times10^{-3}$ & {\rm cm} & Thickness of cathode CL & \\\hline
$d_{\rm Pt}$ & $3\times10^{-7}$ & {\rm cm} & diameter of Pt  nanoparticle & \\\hline
$p_{\rm Pt}$ & $4\times10^{-4}$ & {\rm g/cm$^2$} & Pt loading & \\\hline
$\rho_{\rm Pt}$ & 21.45 & {\rm g/cm$^3$} & density of Pt nanoparticle & \\\hline
$\varepsilon$ & 0.2 & & Volume fraction of ionomer increment in cathode & \cite{LMGAW/15}\\\hline
\multicolumn{5}{|l|}{Physical constants}\\\hline
$R$ & 8.31445985 & {\rm J/mol/K} & Gas constant & \cite{Rum/19}\\\hline
$F$ & 96485.3329 & {\rm C/mol} & Faraday constant & \cite{Rum/19}\\\hline
\multicolumn{5}{|l|}{Parameters for $Pt^{2+}$ formation and diffusion}\\\hline
$\nu_1$ & $1\times10^4$ & {\rm Hz} & Dissolution attempt frequency & \cite{LMGAW/15}\\\hline
$\nu_2$ & $8\times10^5$ & {\rm Hz} & Backward dissolution rate factor & \cite{LMGAW/15}\\\hline
$\beta_1$ & 0.5 & & Butler--Volmer transfer coefficient for Pt dissolution & \cite{LMGAW/15}\\\hline
$n$ & 2 & & Electrons transferred during Pt dissolution &\\\hline
$U_{\rm eq}$ & 1.18 & {\rm V} & Pt dissolution bulk equilibrium voltage & \cite{Dob/75}\\\hline
$\mathit{\Omega}$ & 9.09 & {\rm cm$^3$/mol} & Molar volume of Pt & \cite{LMGAW/15}\\\hline
$\gamma$ & $2.4\times10^{-4}$ & {\rm J/cm$^2$} & Pt [1 1 1] surface tension & \cite{LMGAW/15}\\\hline
$c_{\rm ref}$ & $1$ & {\rm mol/cm$^3$} & reference $Pt^{2+}$ concentration & \cite{HM/12}\\\hline
$H_{1,{\rm fit}}$ & $4\times10^4$ & {\rm J/mol} & Fit Pt dissolution activation enthalpy &  \cite{LMGAW/15}\\\hline
$D_{\rm Pt}$ & $1\times10^{-6}$ & {\rm cm$^2$/s} & Diffusion coefficient of Pt$^{2+}$ in the membrane & \cite{BGADKE/11}\\\hline
\multicolumn{5}{|l|}{Parameters for Pt oxide formation}\\\hline
$pH$ & 0 & & Potential of hydrogen ions (protons) & \\\hline
$\nu_1^\star$ & $1\times10^4$ & {\rm Hz} & Forward Pt oxide formation rate constant &  \cite{LMGAW/15}\\\hline
$\nu_2^\star$ & $2\times10^{-2}$ & {\rm Hz} & Backward Pt oxide formation rate constant &  \cite{LMGAW/15}\\\hline
$\mathit{\Gamma}$ & $2.2\times10^{-9}$ & {\rm mol/cm$^2$} & Pt surface site density &  \cite{LMGAW/15}\\\hline
$\beta_2$ & 0.5 & & Butler--Volmer transfer coefficient for PtO formation & \cite{LMGAW/15}\\\hline
$n_2$ & 2 & & Electrons transferred during Pt oxide formation &\\\hline
$U_{\rm fit}$ & 0.8 & {\rm V} & Pt oxide formation bulk equilibrium voltage & \cite{Dob/75}\\\hline
$\lambda$ & $2\times10^4$ & {\rm J/mol} & Pt oxide dependent kinetic barrier constant &  \cite{LMGAW/15}\\\hline
$\omega$ & $5\times10^4$ & {\rm J/mol} & Pt oxide-oxide interaction energy & \cite{LMGAW/15}\\\hline
$H_{2,{\rm fit}}$ & $1.2\times10^4$ & {\rm J/mol} & Fit partial molar oxide formation activation enthalpy & \cite{LMGAW/15}\\\hline
\end{tabular}
\vspace*{1ex}
\caption{Physical and model parameters applied in the simulation}\label{tab1}
\end{center}
}
\end{table}

\paragraph{Operation of Pt catalyst layer} 

\begin{figure}[hbt!]
\hspace*{-1.5cm}
\epsfig{file=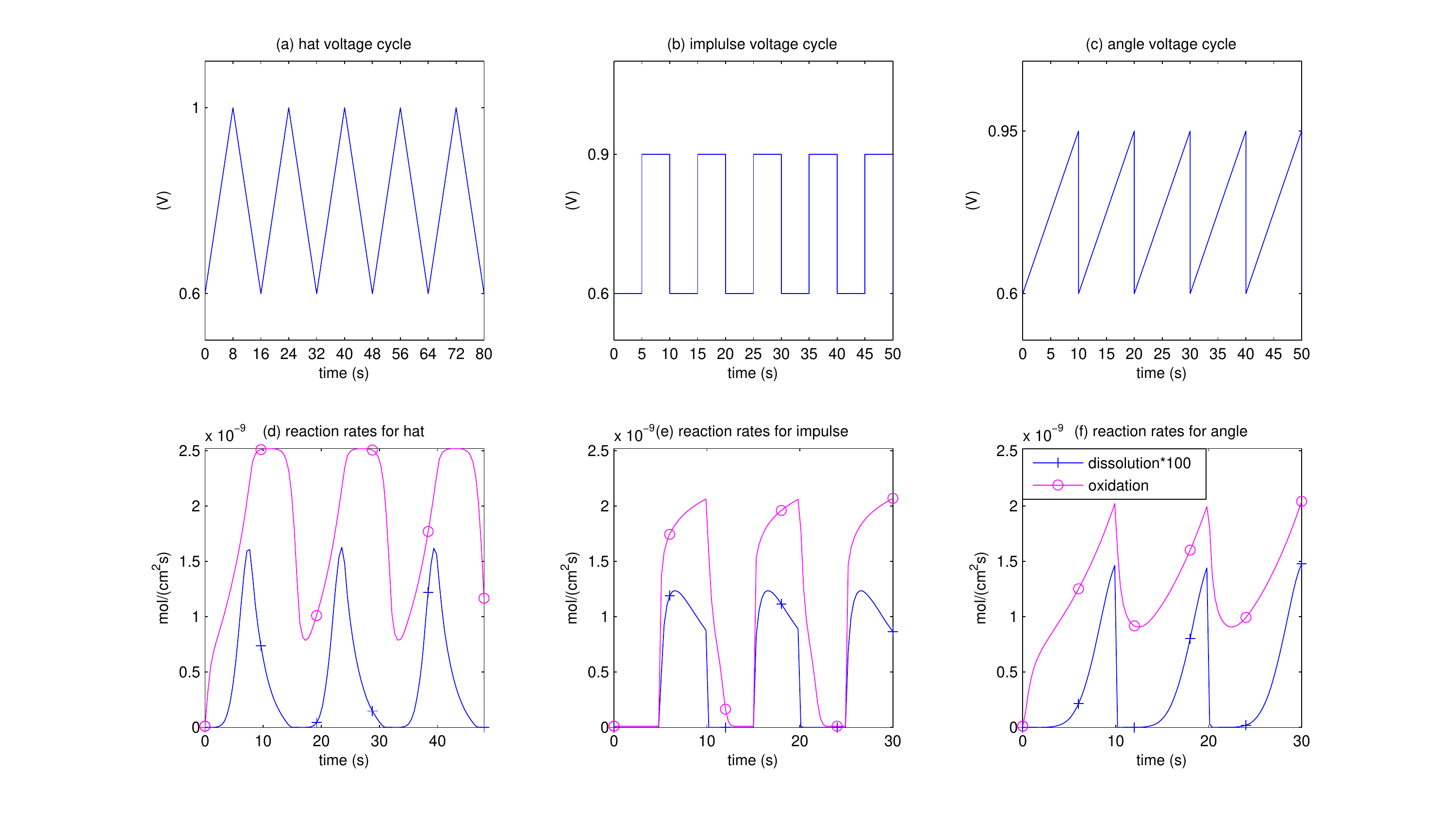,width=1.2\textwidth,angle=0}
\caption{The $\Lambda$-shaped (a) and $\Pi$-shaped (b) 
$\angle$-shaped (c) profiles of cyclic voltage $V(t)$. 
The reaction rates $r_{\rm dissol}$ and 
$r_{\rm oxide}$ for the $\Lambda$ (d), $\Pi$ (e), $\angle$ (f) 
shaped $V(t)$ within 3 cycles.}
\label{fig_voltage_rate_total}
\end{figure}

For our investigation we chose three different protocols 
used by different Institutions to test durability of the catalysts 
applied in low temperature fuel cells: 
\begin{itemize}
\item[(a)]
Accelerated Stress Test used DOE (The U.S. Department of Energy) 
--- $\Lambda$-shaped symmetric triangle wave (50 mV/sec) from 0.6 to 1.0 V 
(see \cite{YBBCGBM/17}); 
\item[(b)]
Accelerated Durability Protocol employed by Tennessee Tech University 
--- $\Pi$-shaped square wave  from 0.6 to 0.9 V, 5 sec at 0.6 V and 5 sec at 0.9 V 
(see \cite{UKRPHLRER/15}); 
\item[(c)]
Durability protocol developed by Nissan --- slow anodic wave 
--- $\angle$-shaped asymmetric triangular wave from 0.6 to 0.95 V 
(see \cite{SMSLMMB/18}). 
\end{itemize}
These profiles are illustrated within 5 periodic cycles in Figure~\ref{fig_voltage_rate_total} 
(a), (b), and (c), respectively.
Each cycle is characterized by the length $p$. 
The $\Lambda$-shaped voltage profile is continuous, symmetric, 
starting and finishing with the minimal voltage value $V_{\rm min}$, 
attaining the maximal voltage value $V_{\rm max}$ 
at the half-length $\frac{p}{2}$, 
thus having the slope $\alpha =\pm\frac{2(V_{\rm max} -V_{\rm min})}{p}$. 
In the plot (a) in Figure~\ref{fig_voltage_rate_total}, $p =16$ (s), 
$V_{\rm min}=0.6$ (V), $V_{\rm max}=1$ (V), $\alpha =\pm5\cdot10^{-2}$ (V/S). 
The $\Pi$-shaped voltage profile is discontinuous, characterized by 
the minimal $V_{\rm min}=0.6$ (V) and the maximal $V_{\rm max}=0.9$ (V) voltages 
switching at $\frac{p}{2}$ given in Figure~\ref{fig_voltage_rate_total} (b) 
for $p =10$ (s). 
The $\angle$-shaped voltage profile first accelerates during $p =10$ (s) 
with the slope $\alpha =3.5\cdot10^{-2}$ (V/S) from the minimal 
$V_{\rm min}=0.6$ (V) to the maximal $V_{\rm max}=0.95$ (V) voltages 
and then switches to $V_{\rm min}$ again, see Figure~\ref{fig_voltage_rate_total} (c). 
\begin{figure}[hbt!]
\begin{center}
\epsfig{file=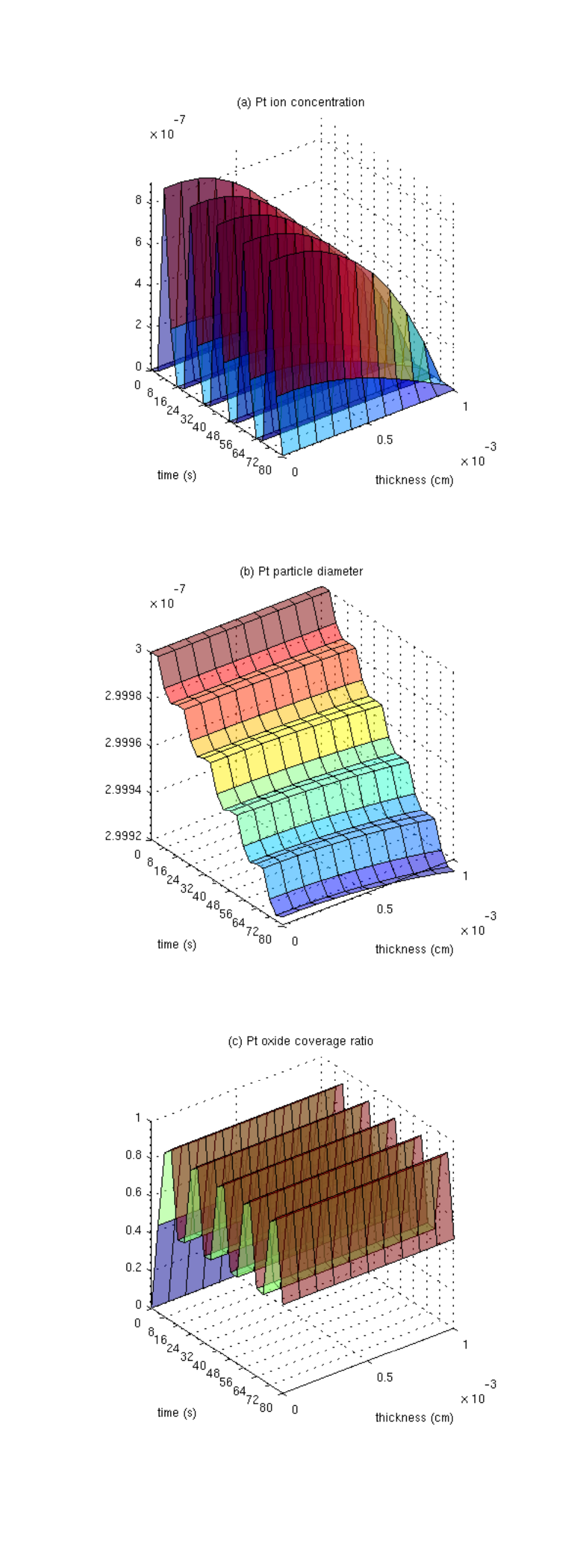,height=.9\textheight,angle=0}
\caption{The solution $c$ in (a), $d$ in (b), $\theta$ in (c) 
under the $\Lambda$-shaped voltage cycle at $T=353.15$ (K).}
\label{fig_cycle}
\end{center}
\end{figure}
\begin{figure}[hbt!]
\begin{center}
\epsfig{file=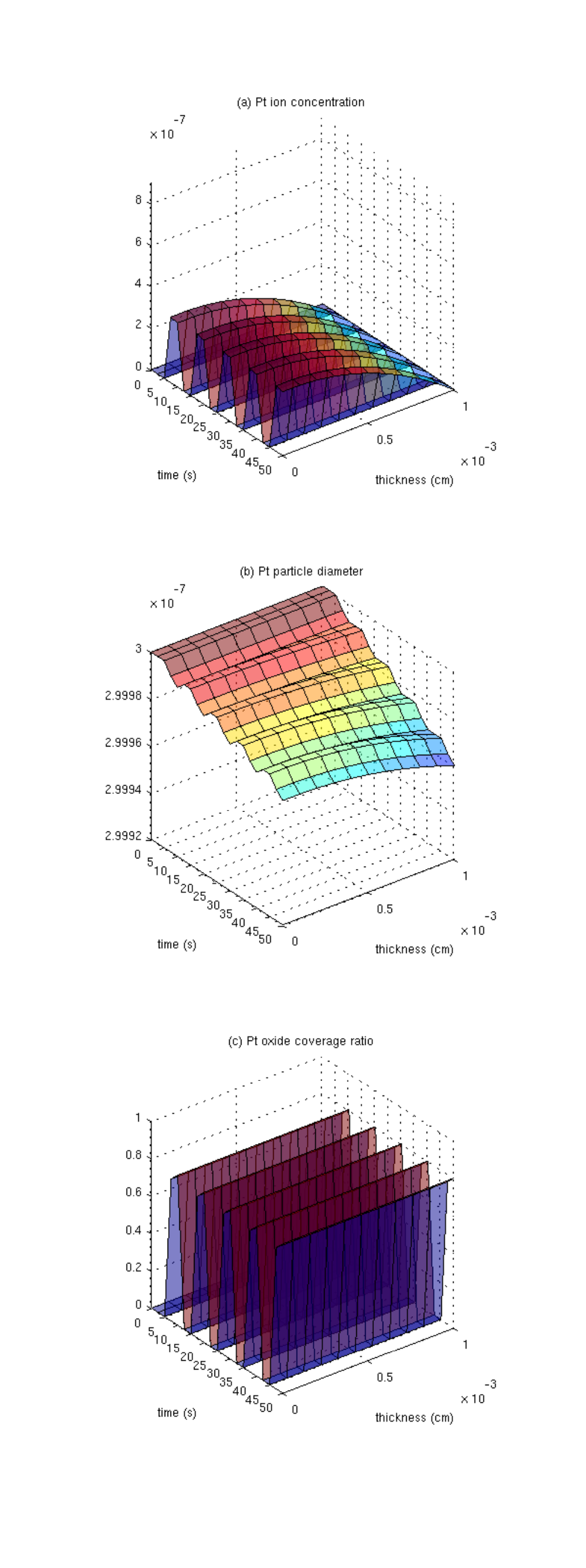,height=.9\textheight,angle=0}
\caption{The solution $c$ in (a), $d$ in (b), $\theta$ in (c) 
under the $\Pi$-shaped voltage cycle at $T=353.15$ (K).}
\label{fig_impulse}
\end{center}
\end{figure}
\begin{figure}[hbt!]
\begin{center}
\epsfig{file=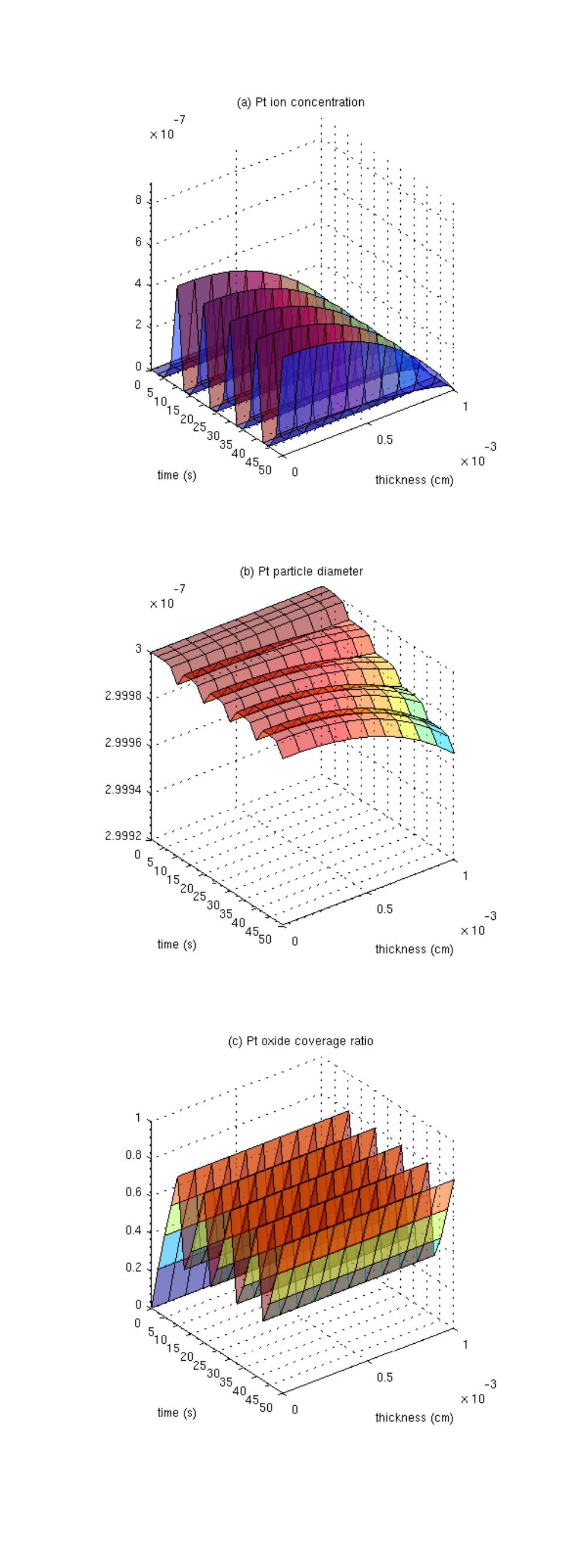,height=.9\textheight,angle=0}
\caption{The solution $c$ in (a), $d$ in (b), $\theta$ in (c) 
under the $\angle$-shaped voltage cycle at $T=353.15$ (K).}
\label{fig_angle}
\end{center}
\end{figure}
In Figure~\ref{fig_voltage_rate_total} (d), (e), (f) we plot the mean over $x\in[0,L]$ in CL of 
the reaction rates $r_{\rm dissol}$ scaled by multiplying by 100, 
and $r_{\rm oxide}$ for the corresponding $\Lambda$, $\Pi$, $\angle$-shaped 
voltage profiles $V(t)$ within 3 periodic cycles.

We use the numerical model \eqref{6} for computer simulation 
of the catalyst coverage and  the CL operation under the $\Lambda$-shaped, 
$\Pi$-shaped, and $\angle$-shaped voltage cycles taken from 
Figure~\ref{fig_voltage_rate_total} (a), (b), and (c). 
The respective solution triples $(c,d,\theta)$ are depicted in 
Figure~\ref{fig_cycle}, \ref{fig_impulse}, \ref{fig_angle} 
at the temperature $T = 80{}^\circ$C. 

From Figures~\ref{fig_cycle}, \ref{fig_impulse}, \ref{fig_angle} (a) 
we observe that at $\Lambda$-shaped voltage cycle the concentration 
of platinum ions is varied from 0 to around $8\times 10^{-7}$ mol/l; 
at $\Pi$-shaped voltage cycle the Pt$^{2+}$ increases maximum to 
$2.3\times 10^{-7}$ mol/l, and at $\angle$-shaped voltage profile it grows 
until $4\times 10^{-7}$ mol/l. 
Figures~\ref{fig_cycle}, \ref{fig_impulse}, \ref{fig_angle} (b) demonstrate 
the evolution in Pt size distribution through the catalyst length. 
The Pt size goes down faster at $\Lambda$-shaped voltage cycle than at other two cycles. 
We should mention, that for all three investigated voltage cycles 
the Pt particle diameter decreases slightly faster at the membrane surface 
than at the boundary with the gas diffusion layer. 
The changes in the coverage of Pt surface by platinum oxide during the voltage cycling 
are depicted in  Figures~\ref{fig_cycle}, \ref{fig_impulse}, \ref{fig_angle} (c). 
In $\Lambda$-shaped and $\angle$-shaped voltage cycles, 
the part of Pt surface is permanently covered by PtO. 
The ratio of the coverage depends on the voltage. 
At $\Lambda$-shaped voltage cycle, the PtO coverage is varied from 42 to 82\%,
while at $\angle$-shaped voltage cycle changes from 30 to 70\%. 
During the $\Pi$-shaped voltage cycle, the platinum oxide covers from 0 to 70\% of Pt surface. 
In this cycle, at the high voltage (0.9 V) the formation of PtO occurs and 
in the next 5 sec, at the low voltage (0.6 V) the reverse reaction proceeds: 
the platinum oxide is reduced to the platinum. 
As seen, at $\Pi$-shaped voltage cycle it is enough time for the reduction reaction, 
while in other two cycles it does not.

\begin{figure}[hbt!]
\begin{center}
\epsfig{file=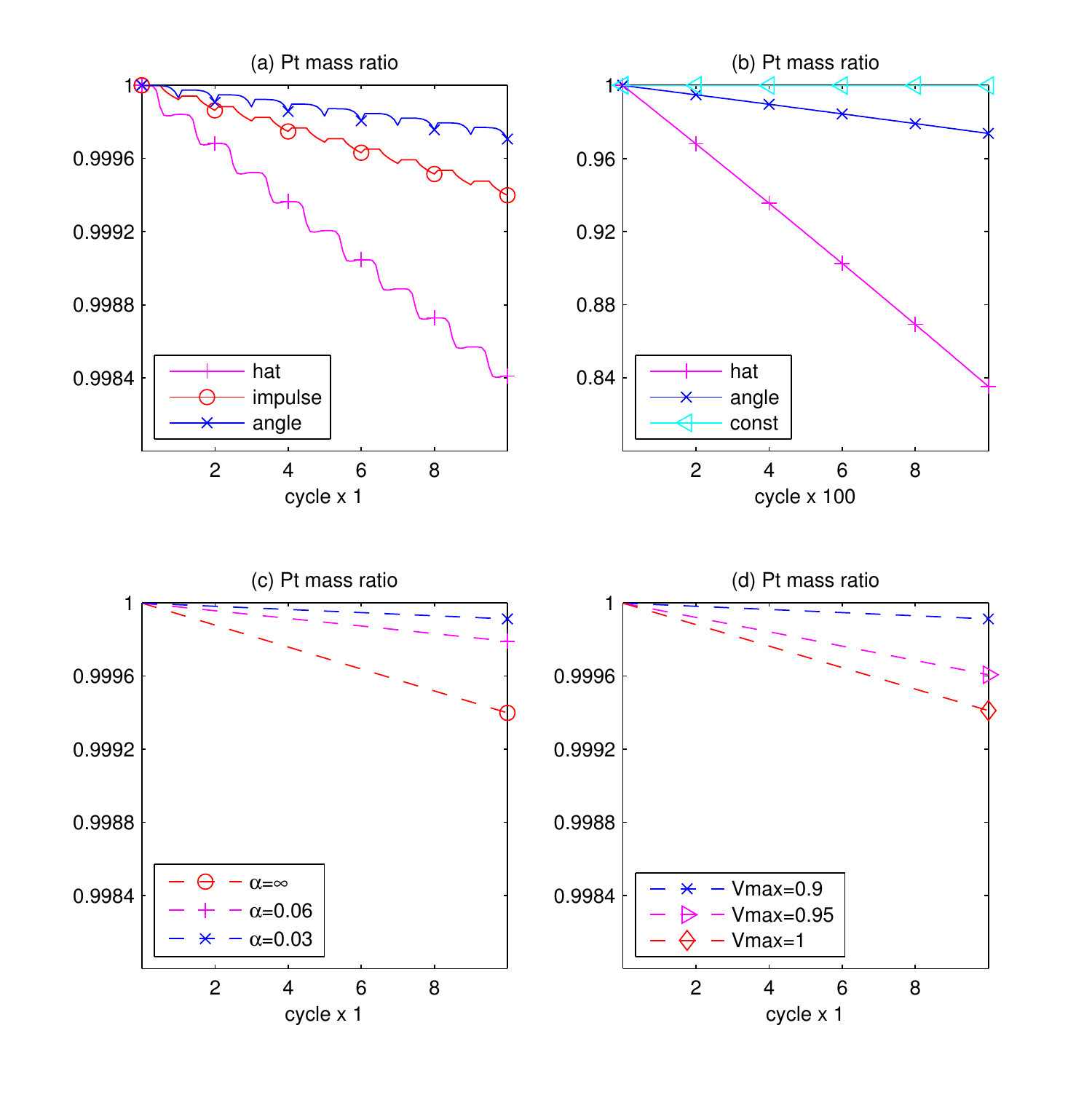,width=.7\textwidth,angle=0}
\caption{The mean Pt mass loss ration $m_{\rm Pt}$ 
under various voltage cycles at $T=353.15$ (K).}
\label{fig_mass_total}
\end{center}
\end{figure}
In Figure~\ref{fig_mass_total} we plot the calculated in time steps $l$ 
platinum mass ratio (the mean over $x\in[0,L]$ in CL): 
\begin{equation}\nonumber
m^l_{\rm Pt} = \frac{4}{3}\pi (d^l)^3/V_{\rm Pt}\in[0,1], 
\end{equation} 
versus the number of cycles. 
For comparison, the three curves $m^l_{\rm Pt}$ are shown 
corresponding to the voltage profiles of hat ($\Lambda$), 
impulse ($\Pi$), and angle ($\angle$) shapes at 10 voltage cycles 
during 2 min. 40 sec. in plot (a). 
While in plot (b) the Pt loss is presented 
at 1000 voltage cycles respectively during 4 hours 26 min. 40 sec. 
There is also depicted $m^l_{\rm Pt}$ under 
the constant voltage $V(t)=0.65$ (V), thus describing idle state. 
The linear loss of Pt mass during voltage cycles can be clear observed 
in Figure~\ref{fig_mass_total}. 
The results shown here are determined by the specific choice 
of profiles and by the upper potential level. 
Indeed, increasing slopes $\alpha =0.03, 0.06, \infty$ (V/s) were tested 
under the fixed upper potential level $V_{\rm max}=0.9$ (V)
as presented in plot (c). 
This confirms that $\Pi$-shaped profile (marked by $\alpha =\infty$) 
is the most damaging with respect to livetimes 
(see experimental data in \cite{UK/07} and \cite{KWSSG/18}, Fig.~4). 
On the other side, increasing the upper potential level 
$V_{\rm max} =0.9, 0.95, 1$ (V) under the fixed slope $\alpha =0.03$ (V/s)
also shortens the lifetime as shown in plot (d) 
(see the experimental confirmation in \cite{KW/19}, Fig.~5(a)). 

\begin{figure}[hbt!]
\hspace*{-1.2cm}
\epsfig{file=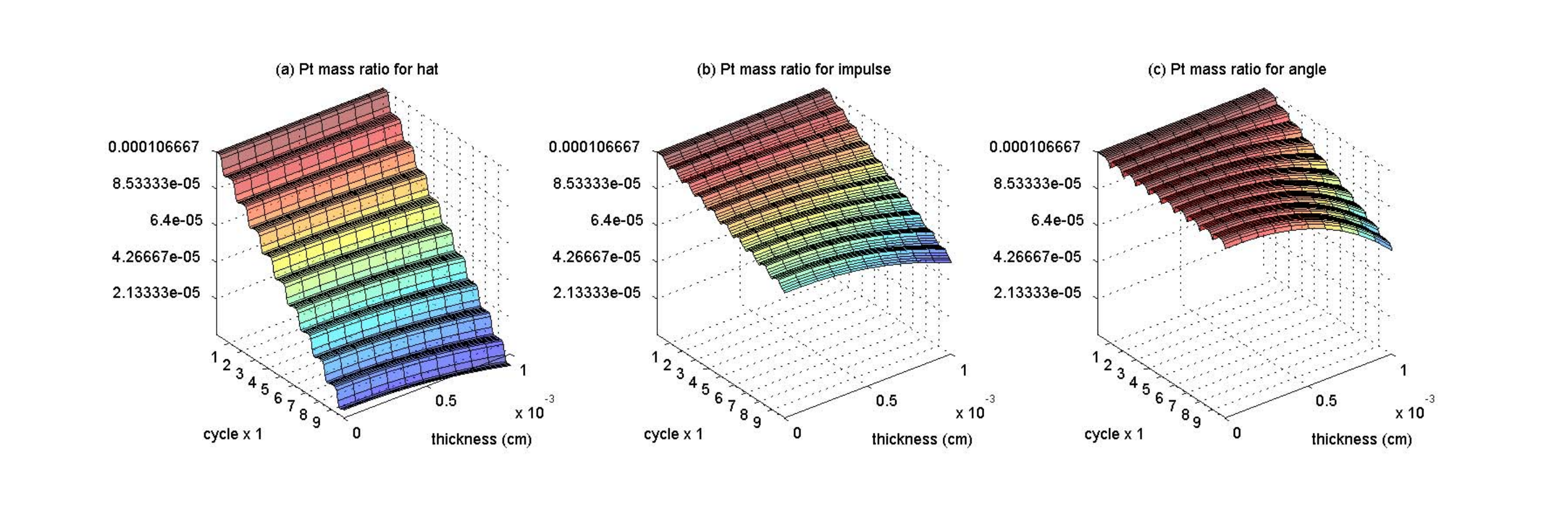,width=1.2\textwidth,angle=0}
\caption{The Pt mass loss ration $m_{\rm Pt}$ vs. CL under $\Lambda$, $\Pi$, 
$\angle$-shaped voltage cycles at $T=353.15$ (K).}
\label{fig_mass_length}
\end{figure}

The degradation phenomenon would be impossible 
without the diffusion (when setting $D_{\rm Pt}=0$ as in \ref{A}). 
In Figure~\ref{fig_mass_length} the corresponding Pt mass loss 
$m^l_{\rm Pt}(x)$ versus $x\in[0,L]$ along CL is presented 
in plots (a), (b), (c) during 10 voltage cycles of $\Lambda$, $\Pi$, $\angle$ profiles. 
Here we can also observe the most strong degradation phenomenon 
near the CL-membrane interface at $x=L$ under the Dirichlet boundary condition. 

Finally, extrapolating the linear Pt loss based on Figure~\ref{fig_mass_total}, 
we calculate the prognosis of the failure when $m_{\rm Pt}$ becomes zero, 
which is presented in Table~\ref{tab2}. 
\begin{table}[hbt!]
{\small
\begin{center}
\begin{tabular}{|l|l|l|l|}\hline
Voltage & Pt mass loss slope & \#cycles prognosis & time prognosis \\\hline
$\Lambda$-shaped & $1.6\times10^{-4}$ & $6\times10^3$ & $27$ h\\\hline
$\Pi$-shaped & $6\times10^{-5}$ & $1.6\times10^4$ & $46$ h\\\hline
$\angle$-shaped & $2.6\times10^{-5}$ & $3.8\times10^4$ & $106$ h\\\hline
constant &$6\times10^{-8}$ & $1.7\times10^7$ & $48000$ h\\\hline
\end{tabular}
\vspace*{1ex}
\caption{The Pt mass loss ration $m_{\rm Pt}$ before failure under various 
voltage cycles at $T=353.15$ (K).}\label{tab2}
\end{center}
}
\end{table}
The average voltages in the studied cycles are 0.8 V for $\Lambda$-shaped cycle; 
0.75 V for $\Pi$-shaped and 0.775 V for $\angle$-shaped cycle. 

\section{Discussion}\label{sec5}

At present, the degradation effects of Pt/C catalyst in PEMFCs 
have been extensively studied as experimentally and theoretically. 
Currently, the decrease in the electroactivity of the Pt/C catalyst 
has been related with the following mechanisms: 
1) Pt dissolution and diffusion into the ionomer; 
2) formation of platinum oxides on the Pt particles surface; 
3) Pt particle ripening; 
4) coalescence of Pt particles. 
The changes in the electrochemical activities of the catalyst are usually 
detected using cycling voltammetry or by measuring polarization curves. 

In our study we consider only two first degradation mechanisms of Pt 
and analyze the Pt mass loss and coverage ratio of Pt particle surface by PtO. 
The simulation presented in the paper is difficult to validate 
with experimental data available in scientific literature. 
However, we found some experimental facts confirming the present calculation. 
For example, Takei et all. \cite{TKKTWU/16} investigated Pt degradation of 
carbon supported Pt catalyst in a single fuel cell as a function of the holding 
times of OCV/load (square voltage cycling) at $T = 80{}^\circ$C and RH 100\%. 
Higher operating potential enhanced the Pt oxidation, 
accelerated the Pt particle growth but suppressed the Pt dissolution. 
The formation of Pt oxide protects the Pt particle from dissolution. 

Ferreira et al. \cite{FlOSHMMKG/05} studied a platinum degradation 
for a short-stack of PEMFC operated at high voltages. 
Using transmission electron microscopy (TEM) they analyzed 
a cross-section of MEA cathode samples and determined 
relative weight percentages of platinum particles 
on the carbon support as a function of cathode thickness. 
The analysis showed that the weight percentage of platinum remaining 
on carbon decreases with increasing distance from gas diffusion layer. 
After the FC operation the weight percentages of Pt at the interface PEM/CL 
was twice time lower than at the interface CL/GDL. 
Yu et al. \cite{YBBCGBM/17} detected 80\% depletion of Pt at cathode/membrane 
interface after accelerated stress test which was performed by 
imposing a triangular wave potential cycling from 0.6 V to 1.0 V 
for 30.000 cycles at 50 mV/sec scan rate. 
These finding confirm our simulation results indicating that the decrease 
in Pt weight at the membrane interface same higher than at the interface to GDL. 

\section{Conclusion}\label{sec6}

We suggested an one-dimensional and dynamic model which describes 
degradation phenomena in Pt catalyst of a polymer electrolyte fuel cell such. 
The model considers Pt dissolution, oxidation as well as diffusion of platinum 
ions through ionomer of the catalyst layer into the membrane. 
Also it takes into account an effect of temperature, Pt particle size and 
Gibbs--Thomson’s effect: dependence of a surface potential on nano-particle size 
as well as influence of the surface potential on the potential gradient in the system. 
The developed model is applied to study concentration profile of Pt$^{2+}$ 
through catalyst length, changes in Pt particle size and mass loss, 
the coverage ratio of Pt surface by platinum oxide at three different voltage cycles 
often used in accelerated stress tests: 
$\Lambda$-shaped, $\Pi$-shaped, and $\angle$-shaped voltage profiles. 

For the parameter values from Table~\ref{tab1}, 
we report here on some of our theoretical and numerical findings with respect 
to admissible voltage operating conditions. 

\begin{itemize}
\item[(i)] 
In order to preserve the physical constraints \eqref{2}, the sufficient are 
$V\in[0.6,1]$ (V) for the accelerated $\Lambda$ and $\angle$-shaped voltage cycles, 
and $V\in[0.6,0.9]$ (V) for the impulse $\Pi$-shaped voltage cycle. 
\item[(ii)] 
In Figures~\ref{fig_cycle}, \ref{fig_impulse}, and \ref{fig_angle} we observe 
diffusion of the Pt ion concentration $c$ and the Pt particle diameter $d$, 
whereas a non-diffusive behavior of the Pt coverage ratio $\theta$. 
\item[(iii)] 
The reaction rates are shown in Figure~\ref{fig_mass_length} for fixed variables, 
and in Figure~\ref{fig_voltage_rate_total} during the CL operation. 
\item[(iv)] 
The rate of the loss of Pt total mass during voltage cycles depends on the 
voltage profile and the upper potential level 
(see Figure~\ref{fig_mass_total} and Figure~\ref{fig_mass_length}, 
and its prognosis in Table~\ref{tab2}). 
\end{itemize}

The study shows that the degradation rate increases with temperature 
and decreasing particle diameter of Pt nano-particles. 
The mass loss in platinum and decrease in Pt particle diameter are 
more significant at the membrane surface than at gas diffusion layer.  

\paragraph{Acknowledgments}
L. K.-J. is supported by the Austrian Research Promotion Agency (FFG) and 
the Austrian Ministry for Transport, Innovation and Technology (BMVIT).\\
V. A. K. is supported by the Austrian Science Fund (FWF) project P26147-N26: PION 
and the European Research Council (ERC) under European Union's Horizon 2020 
Research and Innovation Programme (advanced grant No. 668998 OCLOC),  
he thanks the Russian Foundation for Basic Research (RFBR) 
project 18-29-10007 for partial support. 


\bibliography{kjkbibfile}

\appendix

\section{Non-diffusive case}\label{A}

The non-diffusive case of \eqref{5} is described by the following 
system of nonlinear reaction equations: 
find a triple $c(t)\ge0$, $d(t)\ge0$, $0\le\theta(t)\le1$ such that 
\begin{subequations}\label{A1}
\begin{equation}\label{A1a}
{\textstyle\frac{{\rm d} c}{{\rm d} t}} =B_3 d^2 \bigl( 
B_1(d, \theta) e^{(1-\beta_1) B_4(d, \theta) V} 
-c B_2(d, \theta) e^{-\beta_1 B_4(d, \theta) V} \bigr)
\quad \text{for $t>0$},
\end{equation}
\begin{equation}\label{A1b}
{\textstyle\frac{{\rm d} d}{{\rm d} t}} =-\mathit{\Omega} 
\bigl( B_1(d, \theta) e^{(1-\beta_1) B_4(d, \theta) V} 
-c B_2(d, \theta) e^{-\beta_1 B_4(d, \theta) V} \bigr)
\quad \text{for $t>0$},
\end{equation}
\begin{equation}\label{A1c}
{\textstyle\frac{{\rm d} \theta}{{\rm d} t}} 
+{\textstyle\frac{2 \theta}{d}} {\textstyle\frac{{\rm d} d}{{\rm d} t}} 
={\textstyle\frac{r_{\rm oxide}(\theta)}{\mathit{\Gamma}}}
\quad \text{for $t>0$},
\end{equation}
which is endowed with the initial conditions:  
\begin{equation}\label{A1d}
c(0) =0,\quad d(0) =d_{\rm Pt},\quad \theta(0) =0.
\end{equation}
\end{subequations}
In equations \eqref{A1a} and \eqref{A1b}, 
the expression \eqref{3a} was inserted in $r_{\rm dissol}$. 

Multiplying \eqref{A1a} with $\mathit{\Omega}$ and \eqref{A1b} 
with $B_3 d^2$, after summation we obtain the homogeneous 
ordinary differential equations (ODE): 
\begin{equation}\nonumber
{\textstyle\frac{{\rm d}}{{\rm d} t}} \bigl( 
\mathit{\Omega} c +{\textstyle\frac{B_3}{3}} d^3 \bigr) =0
\quad \text{for $t>0$},
\end{equation}
which implies the first integral of the system  
\begin{equation}\label{A2}
c(t) =c(0) +{\textstyle\frac{B_3}{3\mathit{\Omega}}} d^3(0) 
-{\textstyle\frac{B_3}{3\mathit{\Omega}}} d^3(t).
\end{equation}
With the help of \eqref{A2}, the system \eqref{A1a}--\eqref{A1c} 
can be reduced to two equations for either $c$ and $\theta$, 
or $d$ and $\theta$. 

We note also the following special cases. 
Dividing \eqref{A1c} with $\theta\not=0$ implies the ODE: 
\begin{equation}\label{A3}
{\textstyle\frac{{\rm d}}{{\rm d} t}} \bigl( \ln ( \theta d^2 ) \bigr) 
={\textstyle\frac{r_{\rm oxide}(\theta)}{\mathit{\Gamma} \theta}}
\quad \text{for $t>0$}.
\end{equation}
If the Pt oxidation reaction rate is $r_{\rm oxide}(\theta) \equiv0$ 
(the equilibrium state), then the equation \eqref{A3} is solved trivially as 
\begin{equation}\nonumber
\theta(t) ={\textstyle\frac{\theta(0)}{d^2(t)}} d^2(0) \equiv0
\end{equation}
due to the initial conditions \eqref{A1d}. 

On the other side, if $\theta(t)\equiv1$ then 
$B_1(d, 1) =B_2(d, 1) =0$ in \eqref{3b}. 
Henceforth, the equations \eqref{A1a} and \eqref{A1b} have zero 
right-hand sides and possess the constant solutions $c(t) \equiv c(0) =0$ 
and $d(t) \equiv d(0) =d_{\rm Pt}$. 

\section{Analytical solution}\label{B}
We assume that the Pt particle diameter $d$ and PtO coverage ratio $\theta$ 
are constant in time, hence the coefficients $B_1$, $B_2$, and $B_4$ are constant, 
and \eqref{A1} is reduced to the following Cauchy problem 
for a single linear inhomogeneous ODE: 
\begin{subequations}\label{B1}
\begin{equation}\label{B1a}
{\textstyle\frac{{\rm d} c}{{\rm d} t}}(t)
=B_1 B_3 d^2  e^{(1-\beta_1) B_4 V(t)} 
-c(t) B_2 B_3 d^2  e^{-\beta_1 B_4 V(t)}\quad \text{for $t>0$},
\end{equation}
\begin{equation}\label{B1b}
c(0) =0.
\end{equation}
\end{subequations}
In the following we solve \eqref{B1} analytically 
for non-steady state voltage $V(t)$. 

We introduce an auxiliary function $K(t)$ such that \eqref{B1a} 
can be rewritten equivalently as the following system for $t>0$: 
\begin{subequations}\label{B2}
\begin{equation}\label{B2a}
{\textstyle\frac{{\rm d} K}{{\rm d} t}}(t) 
=B_2 B_3 d^2  K(t) e^{-\beta_1 B_4 V(t)},
\end{equation}
\begin{equation}\label{B2b}
{\textstyle\frac{{\rm d} (K c)}{{\rm d} t}}(t)
=K {\textstyle\frac{{\rm d} c}{{\rm d} t}}(t) 
+c {\textstyle\frac{{\rm d} K}{{\rm d} t}}(t)
=B_1 B_3 d^2 K(t) e^{(1-\beta_1) B_4 V(t)}.
\end{equation}
\end{subequations}
The integration of the equations \eqref{B2a} and \eqref{B2b}, respectively, 
\begin{subequations}\label{B3}
\begin{equation}\label{B3a}
K(t) = K(0) e^{B_2 B_3 d^2 \int_0^t e^{-\beta_1 B_4 V(s)} \,{\rm d}s},
\end{equation}
\begin{equation}\label{B3b}
(K c)(t) =(K c)(0) +B_1 B_3 d^2 \int_0^t e^{(1-\beta_1) B_4 V(\tau)} 
K(\tau) \,{\rm d}\tau,
\end{equation}
\end{subequations}
and the subsequent substitution of \eqref{B3a} into \eqref{B3b} 
leads to the explicit expression
\begin{multline}\label{B4}
c(t) =e^{-B_2 B_3 d^2  \int_0^t e^{-\beta_1 B_4 V(s)} 
\,{\rm d}s} \bigl( c(0)\\ 
+B_1 B_3 d^2 \int_0^t e^{(1-\beta_1) B_4 V(\tau) 
+B_2 B_3 d^2  \int_0^\tau e^{-\beta_1 B_4 V(s)} 
\,{\rm d}s} \,{\rm d}\tau \bigr),  
\end{multline}
which can be shortened using the homogeneous initial condition \eqref{B1b}. 

\end{document}